\documentclass{elsarticle}

\usepackage{graphicx}
\usepackage{float}
\usepackage{subfigure}
\usepackage{amsfonts}
\usepackage{amsmath}
\usepackage{amsthm}
\usepackage{amssymb}
\usepackage{mathtools}
\usepackage{commath}
\usepackage{bm}
\usepackage{caption}

\usepackage[percent]{overpic}
\usepackage[usenames,dvipsnames]{color}
\def\bigdot{\boldsymbol\cdot}

\usepackage{enumitem}
\usepackage{comment}
\usepackage[T1]{fontenc}
\usepackage[utf8]{inputenc}
\usepackage{multirow}

\usepackage{lineno}
\usepackage{natbib}

	{\left\lbrace \begin{array}{@{} l @{} }}%
	{ \end{array}\right.}

\newtheorem{pro}{Proposition}

\newtheorem{remark}{Remark}

\newcommand{\bs}{\boldsymbol}

\newlist{level}{itemize}{4}
\setlist[level]{label={},noitemsep,topsep=0pt}
\newcounter{algo}
\renewcommand{\thealgo}{\arabic{algo}}
\newenvironment{algorithm}[1]{%
    \refstepcounter{algo}%
    \paragraph{Algorithm \thealgo}#1%
    \vspace{2pt}\hrule\vspace{5pt}%
    \begin{level}
}{%
    \end{level}%
    \vspace{5pt}\hrule\vspace{\baselineskip}%
}

\makeatletter
\def\algbackskip{\hskip-\ALG@thistlm}
\makeatother
\begin{document}

\begin{frontmatter}



\title{A 
stochastic extended Rippa's algorithm for LpOCV}


\author{L. Ling$^{*}$}
\ead{lling@hkbu.edu.hk}

\author{F. Marchetti$^{**}$}
\ead{francesco.marchetti@unipd.it}

\address{$^{*}$Department of Mathematics, Hong Kong Baptist University, Hong Kong; \\$^{**}$Istituto Nazionale di Alta Matematica (INdAM) - Dipartimento di Matematica \lq\lq Tullio Levi-Civita\rq\rq, Universit\`a di Padova, Italy}

\begin{abstract}

In kernel-based approximation, the tuning of the so-called \textit{shape parameter} is a fundamental step for achieving an accurate reconstruction. Recently, the popular Rippa's algorithm \cite{Rippa99} has been extended to a more general cross validation setting. In this work, we propose a modification of such extension with the aim of further reducing the computational costs. The resulting Stochastic Extended Rippa's Algorithm (SERA) is first detailed and then tested by means of various numerical experiments, which show its efficacy and effectiveness in different approximation settings.
\end{abstract}

\begin{keyword}
RBF interpolation, Extended Rippa's Algorithm, cross validation, stochastic low-rank approximation

\MSC[2020]: 65D05, 65D12, 41A05
\end{keyword}

\end{frontmatter}

\section{Introduction}\label{sec:intro}

Let $\Omega\subset \mathbb{R}^d$ 
be some bounded domain with sufficiently smooth boundary. For some  \textit{shape parameter} $\varepsilon>0$, we consider a strictly positive definite, radial, and translation invariant kernel $\kappa_{\varepsilon}:\Omega\times\Omega\longrightarrow\mathbb{R}$
in the form of
$$\kappa_{\varepsilon}(\bs{x},\bs{y})=\varphi_{\varepsilon}(\lVert \bs{x}-\bs{y} \lVert )  
\qquad \mbox{ for }\bs{x},\bs{y}\in\Omega,$$
where
$\lVert\bigdot\lVert:= \lVert\bigdot\lVert_{\ell^2(\mathbb{R}^d)}$ is the Euclidean norm with some  univariate \emph{radial basis function} $\varphi_{\varepsilon}:
[0,\infty) 
\longrightarrow \mathbb{R}$, say, Gaussian RBF $\varphi_{\varepsilon}(r)=\exp(-(\varepsilon r)^2 )$.

We aim to recover an unknown function $f:\Omega\longrightarrow\mathbb{R}$ based on $n\in\mathbb{N}$ nodal values at $X=\{\bs{x}_1,\dots,\bs{x}_n\}$ denoted by $\bs{f}_X=(f(\bs{x}_1),\dots,f(\bs{x}_n))^{\intercal}$ with an ansatz
\begin{equation}\label{eq:interpolant}
    S(\bs{x})=\sum_{i=1}^n c_i\kappa_{\varepsilon}(\bs{x},\bs{x}_i),
\end{equation}
by imposing interpolation conditions at $X$.
The vector of coefficients $\bs{c}=(c_1,\dots,c_n)^{\intercal}$ in (\ref{eq:interpolant}) is uniquely determined as the solution of the linear system
\begin{equation}\label{eq:problem}
    \mathsf{K}_{\varepsilon}(X,X) \bs{c} = \bs{f}_X.
\end{equation}
The entries of the \textit{kernel matrix} in (\ref{eq:problem}) are given by
$[\mathsf{K}_{\varepsilon}(X,X)]_{ij}=\kappa_{\varepsilon}(\bs{x}_i,\bs{x}_j)$ for $\bs{x}_i,\bs{x}_j\in X$. If $\mathsf{K}_{\varepsilon}$ is severely ill-conditioned and some regularization is needed (e.g. Tikhonov regularization), then the interpolation conditions might be relaxed for an improved stability; for a complete introduction concerning kernel-based approximation, refer to e.g. \cite{Fasshauer07,Fasshauer15}.

The quality of a constructed approximant is often strongly influenced by the value of the shape parameter, and consequently many strategies have been proposed in the literature for its tuning; see e.g. \cite{Cavoretto19,Cavoretto18,Fornberg07,Scheuerer11} and \cite[Chapter 14]{Fasshauer15} for an overview.
To perform such tuning, a possibility consists in performing a Cross Validation (CV) scheme \cite{Golub79}. First, the dataset
(i.e., $X$ and $\bs{f}_X$ in our context)
is divided into $k\in\mathbb{N}$ (possibly equal-sized) disjoint subsets, $k\le n$. Then, $k$ different interpolants, a.k.a. models, as in (\ref{eq:interpolant})
are built upon $k-1$ \textit{training} folds.
Their performance is then assessed on the respective remaining \textit{validation} fold. By setting $p\approx n/k$, $p\in\mathbb{N}$, this procedure can be intended as an approximation of \textit{Leave-$p$-Out} CV (L$p$OCV), where all possible combinations of $p$ elements are taken into account as validation set \cite{Celisse08}. In the particular case $p=1$, the resulting \textit{Leave-One-Out} CV (LOOCV) scheme computes an exact $n$-fold CV, and it has been widely employed by the scientific community and also generalized to other contexts e.g. in \cite{Fasshauer07a,Mongillo11}. In this paper, with an abuse of definition, we will refer to L$p$OCV meaning $k$-fold CV with $k\approx n/p$.

The purpose of this work is providing a modification of the Extended Rippa's Algorithm (ERA), which is described in Section \ref{sec:era}, for achieving an even better trade-off between computational time and accuracy of the resulting interpolant. The proposed approach, which is valid also in the more general setting of other single-parameter family of kernel matrix systems, say, Tikhonov regularization, is detailed in Section \ref{sec:sera} and tested in Section \ref{sec:numerics}.

\section{Extended Rippa's Algorithm (ERA)}\label{sec:era}

Assuming $p\ll n$, a \textit{naive} implementation of L$p$OCV leads to a computational cost of order $\mathcal{O}(n^4/p)$ for inverting approximately $k\approx n/p$ different $(n-p)\times(n-p)$ linear systems, i.e., the cost of constructing \textit{circa} $k$ interpolants whose related kernel matrix is of order $n-p$. Recently in \cite{Marchetti21}, an ERA has been proposed for a more efficient computation of L$p$OCV.
The algorithm involves the computation of the matrix inverse $\mathsf{K}_{\varepsilon}^{-1}:=[\mathsf{K}_{\varepsilon}(X,X)]^{-1}$, which is $\mathcal{O}(n^3)$, and the resolution of approximately $n/p$ different $p\times p$ linear systems, which is $\mathcal{O}(np^2)$ (see \cite[Section 2.2]{Marchetti21}).
More precisely, by using the notation of \cite{Marchetti21}, let a value of $\varepsilon$ be fixed; let $\bs{v} =(v_1,\dots,v_p)^{\intercal} \in\{1,\dots,n\}^p$, 
be one of the $k$ vectors of distinct validation indices, then the vector of $\varepsilon$-dependent  \emph{validation errors} $\bs{e}_{\bs{v}}=\bs{e}_{\bs{v}}(\varepsilon)$ related to the nodes $\bs{x}_{v_1},\dots,\bs{x}_{v_p}$ satisfies $[\mathsf{K}_{\varepsilon}^{-1}]_{\bs{v},\bs{v}}\bs{e}_{\bs{v}}=\bs{c}_{\bs{v}}$, where the vectorized index $\bs{v}$ extracts the $i$-th rows and $j$-th columns subsystem with $i,j\in \bs{v}$ of the original.
Augmenting all $k$ validation error  vectors $\bs{e}(\varepsilon):=[\bs{e}_{\bs{v}_1}^{\intercal},\dots,\bs{e}_{\bs{v}_{k}}^{\intercal}]^{\intercal}(\varepsilon)$ to form the full error,
we define the LpOCV optimal value by $\varepsilon^*= \arg\min \lVert \bs{e}(\varepsilon) \lVert$.

\section{Stochastic Extended Rippa's Algorithm (SERA)}\label{sec:sera}

The computational cost of the ERA is dominated by the calculation of $\mathsf{K}_{\varepsilon}^{-1}$. In this section, we propose a stochastic strategy for approximating $\mathsf{K}_{\varepsilon}^{-1}$ with the aim of further speeding up the calculations. The scheme that we propose is inspired by an idea in \cite[Section 3.1]{Yang18}, where the diagonal of $\mathsf{K}_{\varepsilon}^{-1}$ is estimated by means of a stochastic procedure, in the framework of LOOCV Rippa's scheme (see also \cite{Bekas07}). In our more general L$p$OCV setting, many submatrices of $\mathsf{K}_{\varepsilon}^{-1}$ need to be approximated, thus some modifications are required.

Therefore, we propose the following stochastic low-rank approximation of $\mathsf{K}_{\varepsilon}^{-1}$ (see Algorithm \ref{alg1}), and we refer to the resulting scheme that makes use of such an approximation as SERA.
\begin{algorithm}{Inverse matrix approximation in SERA}\label{alg1}
    \item \textbf{Inputs:}\\
    $\mathsf{K}_{\varepsilon}$: $n\times n$ kernel matrix;\\
    $s$: natural number, $0<s<n$.
    \item \textbf{Core}
    \begin{level}
        \item
        1. Generate an $n\times s$ matrix $\mathsf{W}_s$ whose entries follow normal distribution $w_{ij} \sim\mathcal{N}(0,1)$.\\
        2. Define $\mathsf{U}_s=\mathsf{K}_{\varepsilon}\mathsf{W}_s$.\\
        3. Calculate the Moore–Penrose inverse $\mathsf{U}^{+}_s:=(\mathsf{U}_s^{\intercal}\mathsf{U}_s)^{-1}\mathsf{U}_s^{\intercal}$, say, by SVD.\\
        4. Compute $\mathsf{V}_s=\mathsf{W}_s \mathsf{U}_s^{+}$ and employ it as an approximation of $\mathsf{K}_{\varepsilon}^{-1}$.
    \end{level}
    \item \textbf{Output:}\\
    $\mathsf{V}_s$: $n\times n$ matrix approximating $\mathsf{K}_{\varepsilon}^{-1}$.
\end{algorithm}

\begin{remark}
    The column space of $\mathsf{U}_s$ is a random subspace of the column space of $\mathsf{K}_{\varepsilon}$. Furthermore, we point out that $\mathsf{V}_s=\mathsf{W}_s \mathsf{U}_s^{+}=\mathsf{K}_{\varepsilon}^{-1}\mathsf{U}_s\mathsf{U}_s^{+}$, and $\mathsf{U}_s\mathsf{U}_s^{+}$ is the orthogonal projector onto the column space of $\mathsf{U}_s$.
\end{remark}

\begin{pro}
The matrix $\mathsf{V}_s$ in Algorithm \ref{alg1} has rank at most $s$. Moreover, the computational cost required for its calculation is $\mathcal{O}(s^2n+sn^2)$.
\end{pro}
\begin{proof}
The first claim is ensured by a property concerning the rank of the product of matrices. Then, Step 2 is $\mathcal{O}(n^2s)$, Step 3 is $\mathcal{O}(ns^2)$ (using SVD) and finally Step 4 is $\mathcal{O}(n^2s)$.
\end{proof}

The overall complexity of SERA is then $\mathcal{O}(s^2n+sn^2+np^2)$ (see Section \ref{sec:intro}), which is quadratic with respect to $n$. Therefore, comparing the stochastic approximation to the exact ERA that is $\mathcal{O}(n^3)$, there is convenience if $s$ is relatively small with respect to $n$. However, the smaller the value of $s$, the poorer the approximation quality of $\mathsf{V}_s$, with a possible loss in terms of accuracy. Besides investigating these two intertwined aspects, in the next section we show that the results obtained by the stochastic approach on a fine evaluation (test) grid are comparable to the ones attained by the exact scheme.

\section{Numerics}\label{sec:numerics}

The experiments have been carried out in Matlab on a Intel(R) Core(TM) i7-1165G7 CPU@2.80GHz processor. The Matlab software with the implementation of the proposed scheme is available for the scientific community at
\begin{equation*}
\texttt{https://github.com/cesc14/RippaExtCV}\:.
\end{equation*}
In order to provide statistical significance, for each considered value of $s$, we run SERA $100$ times and we report the results by using boxplots; on each box, the central red mark indicates the median, the bottom and top edges of the box indicate the 25th and 75th percentiles, respectively, the whiskers stretch to the most extreme data points not considered outliers, and finally the outliers are displayed by using red crosses. We observe that if we set $s=\lfloor \sqrt{n}\rfloor$ and we repeat SERA $\lfloor\sqrt{n}\rfloor$ time, then the overall cost is still $\mathcal{O}(n^3)$, i.e., the same order of ERA.

\subsection{ERA vs. SERA}\label{sec:testuno}

Let $\Omega=[-1,1]^2$ and let $f:\Omega\longrightarrow\mathbb{R}$ be the test function
in \cite{Marchetti21} defined as
\begin{equation*}
    f(\bs{x})=\frac{\sin(x_1)}{x_1^2+1}\cdot \frac{\cos(x_2)}{x_2^2+1},\quad \bs{x}=(x_1,x_2).
\end{equation*}
In the following, we test the proposed scheme in two numerical experiments. In both tests, we take the classical $L_2$-norm for the error and we tune the shape parameter $\varepsilon$ by using both ERA and SERA. After that, we test the performance of the approximation schemes resulting from the parameter tuning by evaluating the error with respect to $f$ on a test equispaced $m\times m$ grid $\Xi_m$ in $\Omega$, with $m=30$. Concerning SERA, the entries of $\mathsf{W}_s$ are normal random numbers. Moreover, letting $q_1,\dots,q_{\ell}$ be equispaced \textit{reduction ratio} values between $q_1=0.05$ and $q_{\ell}=0.5$, we set $\ell = 16$ and we let $s$ varying between $s_1=\lfloor 0.05 n^2 \rfloor$ and $s_{\ell}=\lfloor 0.5 n^2 \rfloor$ by considering the values $s_i=\lfloor q_i n^2\rfloor$, $i=1,\dots,\ell$.

\subsubsection{Test 1a: \textit{Matérn $C^0$} kernel}
\begin{figure}[h!]
     \centering
    \subfigure[The chosen $\varepsilon^{\star}\in\bs{\varepsilon}_1$.]{\label{fig:1a}\includegraphics[width=0.327\linewidth]{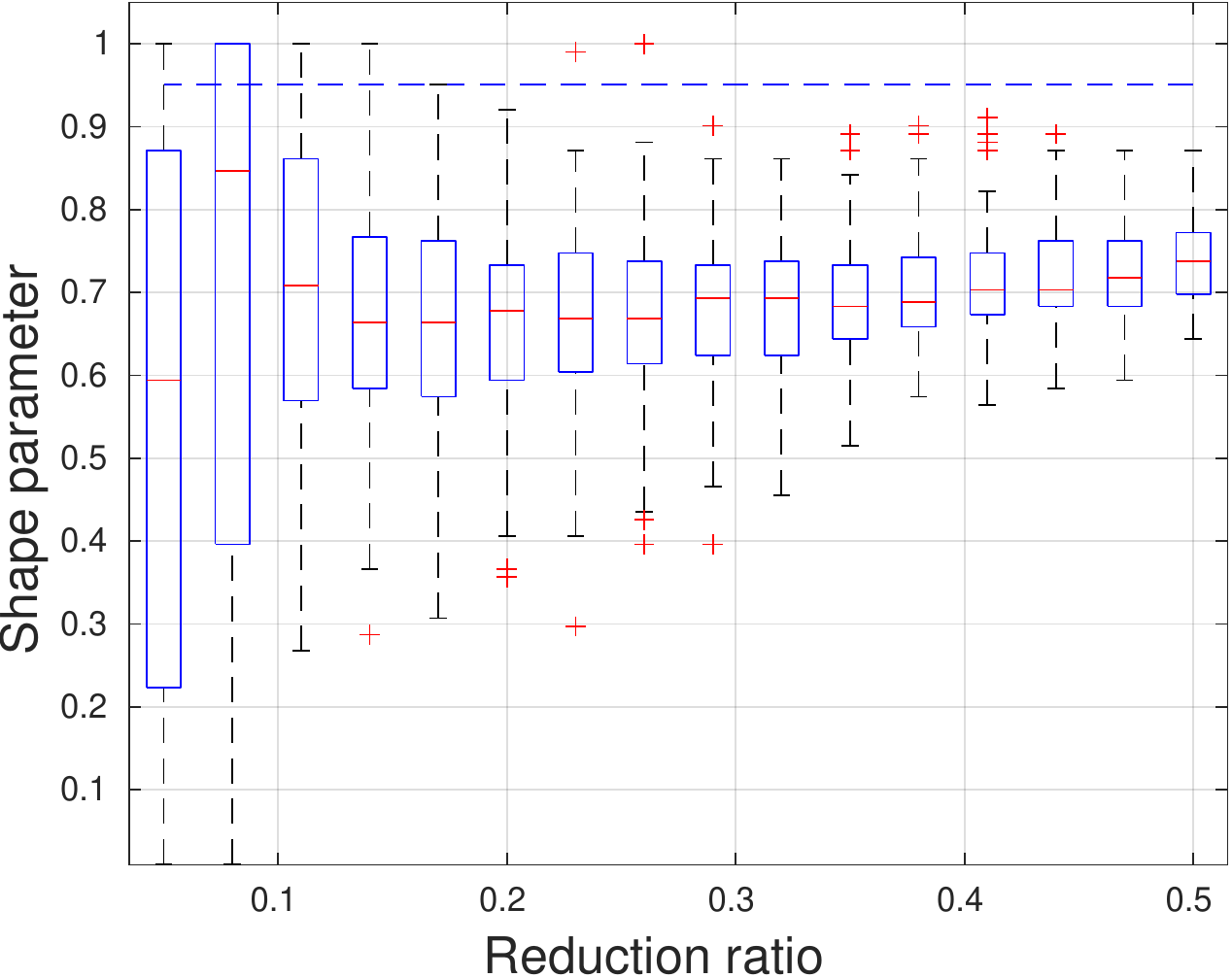}}
    \subfigure[Time employed.]{\label{fig:1b}\includegraphics[width=0.327\linewidth]{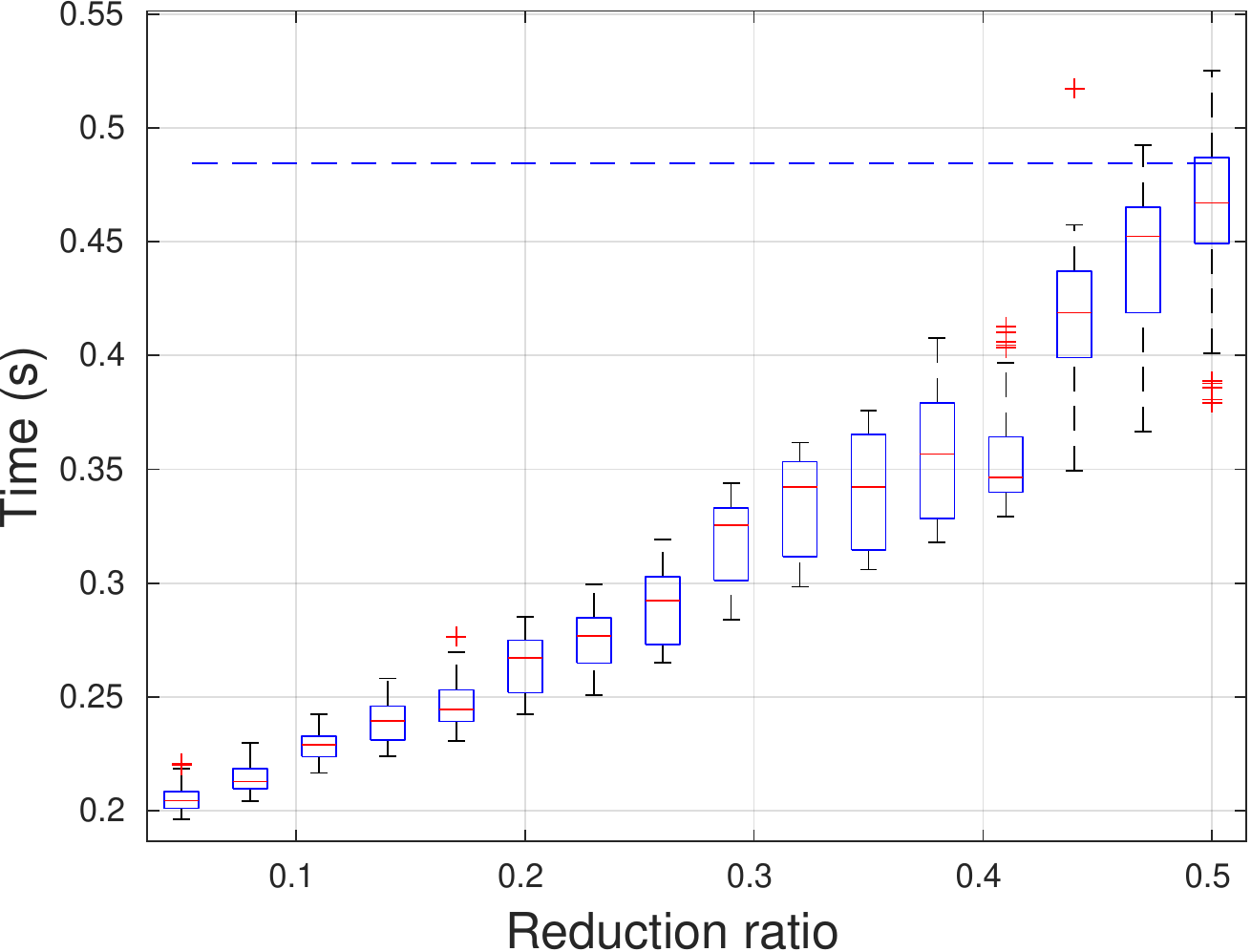}}
    \subfigure[Test error on $\Xi_m$.]{\label{fig:1d}\includegraphics[width=0.327\linewidth]{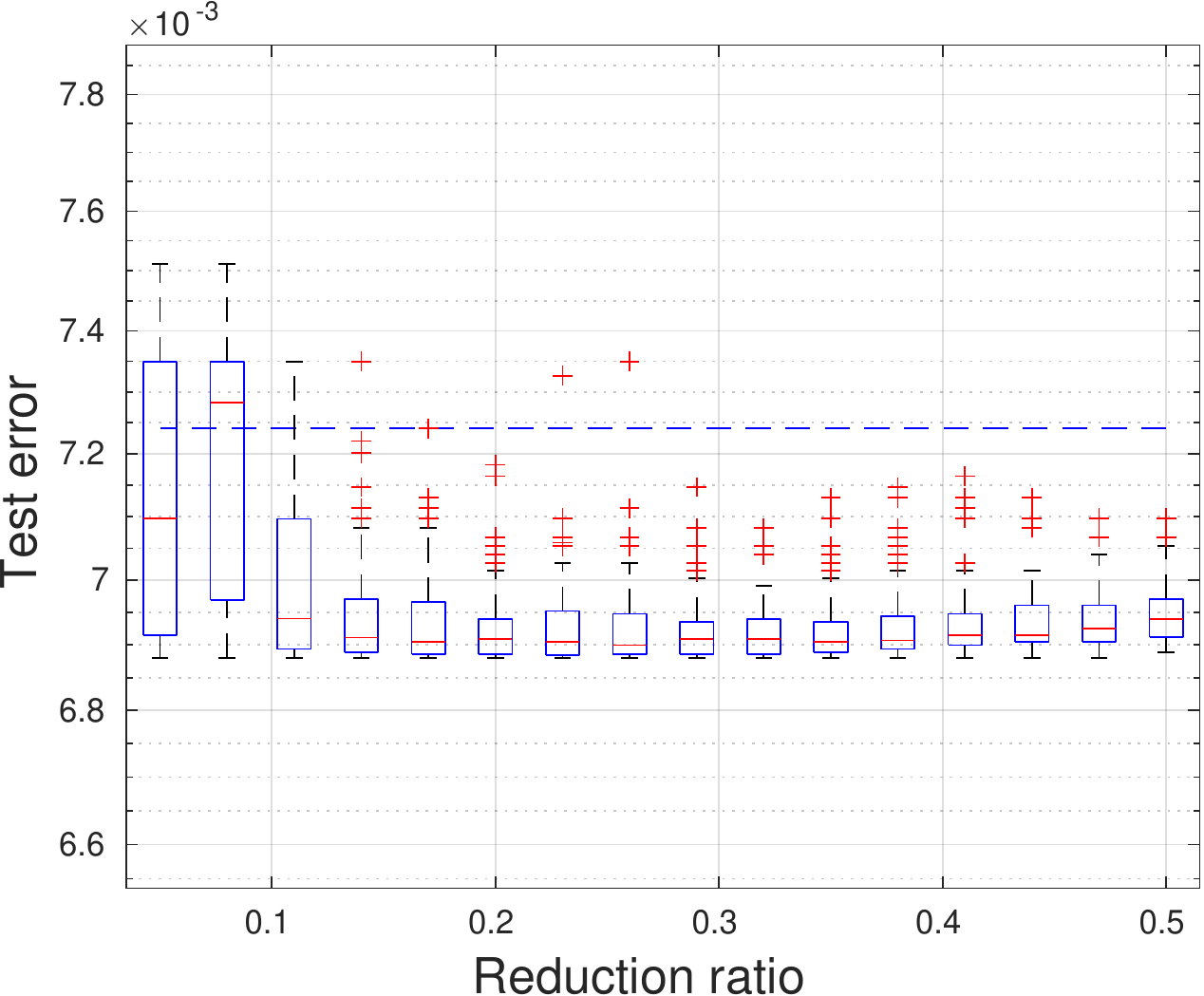}}
        \caption{Test 1a. On the $x$-axis, we let $q_1,\dots,q_{\ell}$ vary. The results related to ERA are depicted by means of a dashed blue line, while the performance of SERA varying $s_i$, $i=1,\dots,\ell$, are reported by using boxplots. We consider the kernel $\varphi_{M,\varepsilon}$ and the interpolation set $E_n$.}
        \label{fig:1}
\end{figure}

Similarly, we consider the following \textit{Matérn $C^0$} kernel (see e.g. \cite[Section 4.4]{Fasshauer07}) $\varphi_{M,\varepsilon}(r)=e^{-\varepsilon r}$, and the equispaced grid $n\times n$ $E_n$ of interpolation nodes in $\Omega$, with $n=20$.
We set $p=2$, i.e. we consider approximated L2OCV, and we take $\bs{\varepsilon}_1 \in [0.01,1]$, uniformly, as vector of $101$ shape parameter values. The results obtained by using ERA and SERA are displayed in Figure~\ref{fig:1}. More precisely, we report the shape parameter $\varepsilon^{\star}$ chosen during validation (Figure \ref{fig:1a}), the time employed by ERA and SERA (Figure \ref{fig:1b}) and the test error with respect to $f$ on $\Xi_m$ achieved by interpolating at $E_n$ with $\varphi_{M,\varepsilon^{\star}}$ (Figure \ref{fig:1d}).

\subsubsection{Test 1b: Gaussian kernel}\label{sec:gaussian_test}
\begin{figure}[h!]
     \centering
    \subfigure[The chosen $\varepsilon^{\star}\in\bs{\varepsilon}_2$.]{\label{fig:2a}\includegraphics[width=0.327\linewidth]{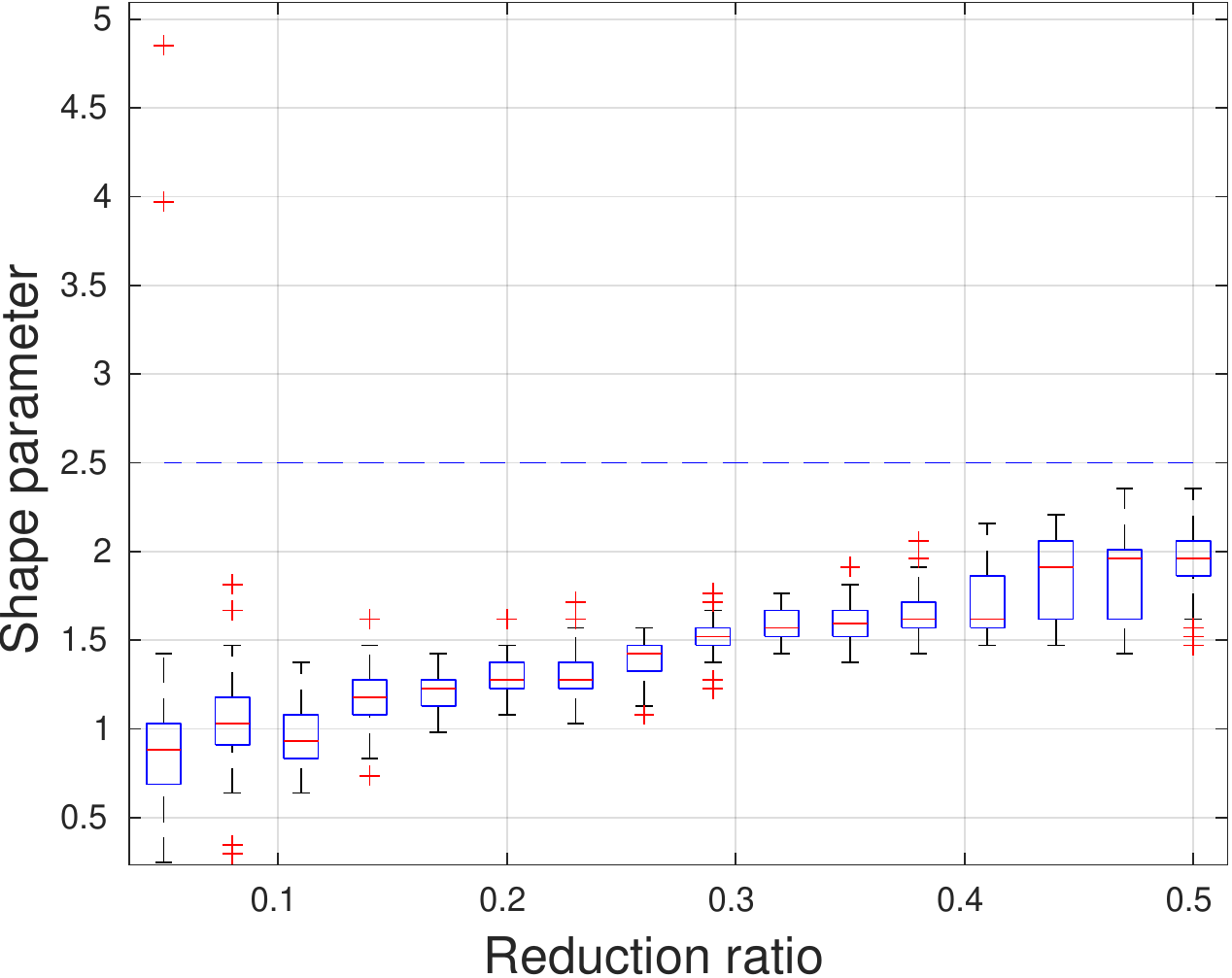}}
    \subfigure[Time employed.]{\label{fig:2b}\includegraphics[width=0.327\linewidth]{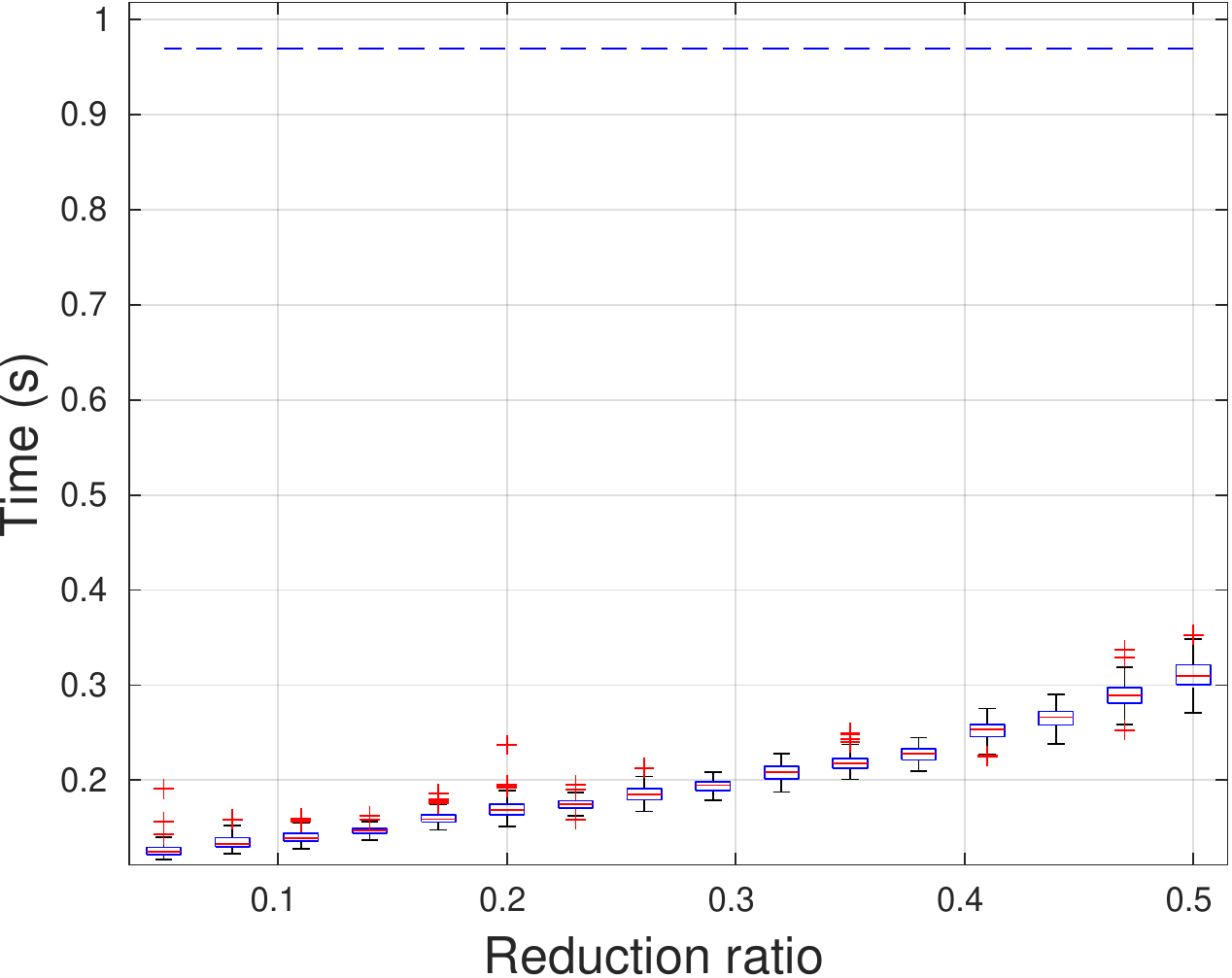}}
    \subfigure[Test error on $\Xi_m$.]{\label{fig:2d}\includegraphics[width=0.327\linewidth]{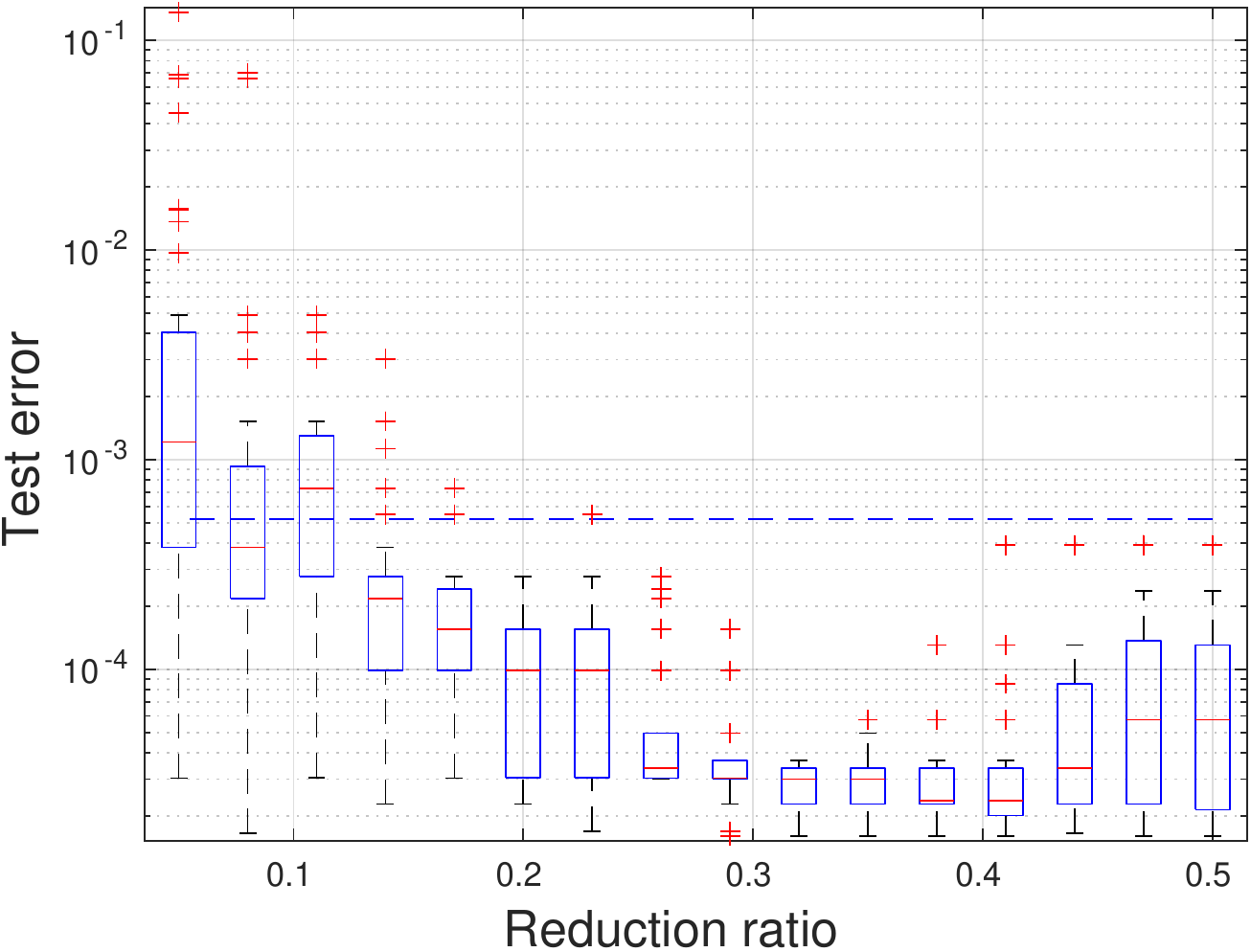}}
        \caption{Test 1b. On the $x$-axis, we let $q_1,\dots,q_{\ell}$ vary. The results related to ERA are depicted by means of a dashed blue line, while the performance of SERA varying $s_i$, $i=1,\dots,\ell$, are reported by using boxplots. We consider the kernel $\varphi_{G,\varepsilon}$ and the interpolation set $H_n$.}
        \label{fig:2}
\end{figure}

We consider the well-known Gaussian kernel (see e.g. \cite[Section 2.1]{Fasshauer07}) $\varphi_{G,\varepsilon}(r)=e^{-(\varepsilon r)^2}$. We set $p=5$ and we take $\bs{\varepsilon}_2 \in [0.1,10]$, uniformly, as vector of $101$ shape parameter values that are considered in the validation process. The interpolation set is the set of quasi-random $n^2$ \textit{Halton points} \cite{Halton60} $H_{n}$ in $\Omega$, with $n=17$. Due to the choice of Gaussian kernel, which often leads to ill-conditioned interpolation processes (see e.g. \cite[Section 2.1]{Fasshauer07}), we have to employ some regularization strategy. Therefore, we use QR factorization in ERA and SERA for regularizing the matrix inversions and the greedy algorithm provided \cite{Ling16} for the construction of the interpolant after the tuning of the shape parameter. Analogously to Figure \ref{fig:1}, in Figure \ref{fig:2} we display the results achieved in this experiment setting. Here, with the considered regularization techniques, the computational advantage with respect to ERA is even more remarkable.

\subsection{Test 2: SERA varying $n$}

\begin{figure}[h!]
     \centering
    \subfigure[The chosen $\varepsilon^{\star}\in\bs{\varepsilon}_3$.]{\label{fig:3a}\includegraphics[width=0.327\linewidth]{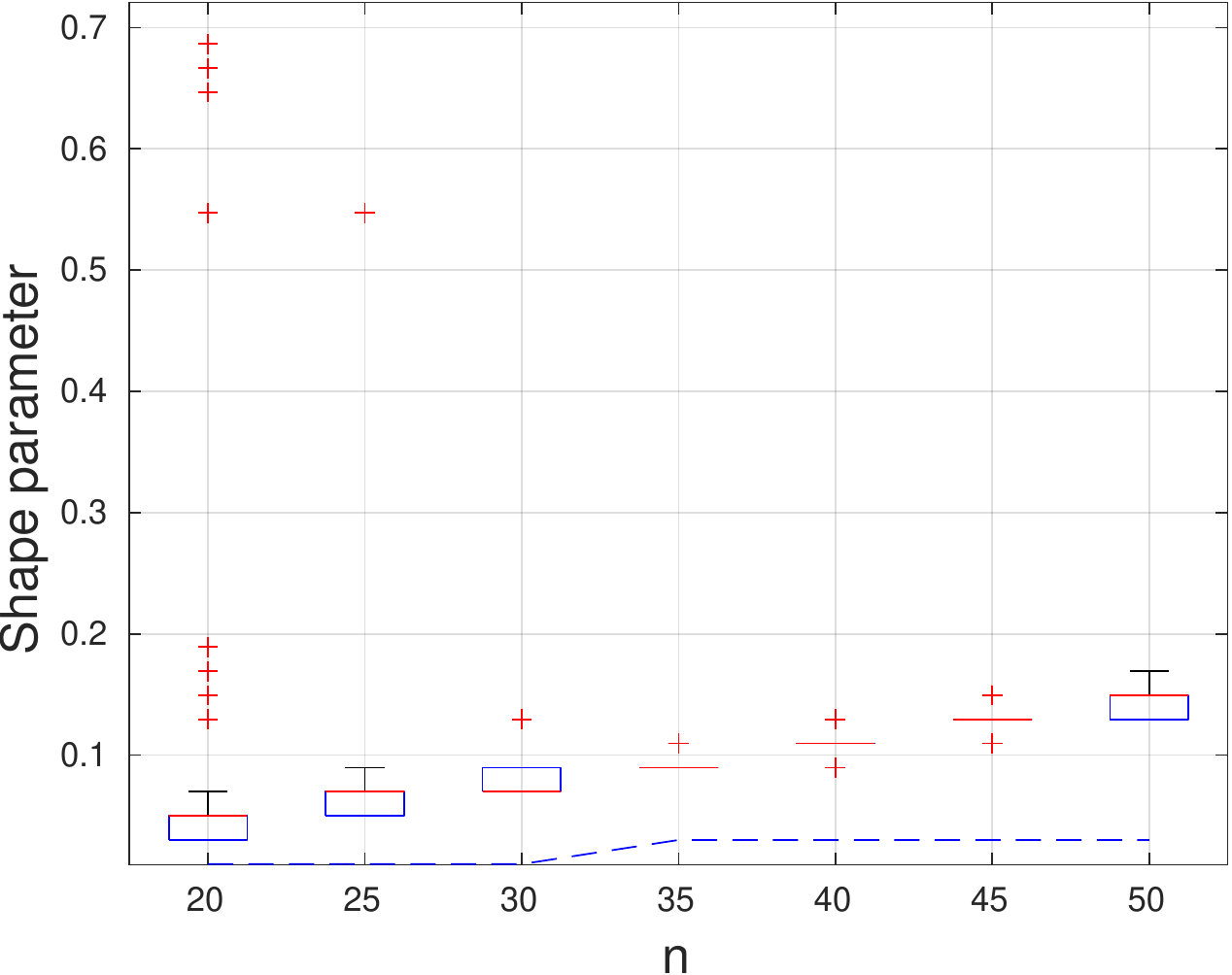}}
    \subfigure[Time employed.]{\label{fig:3b}\includegraphics[width=0.327\linewidth]{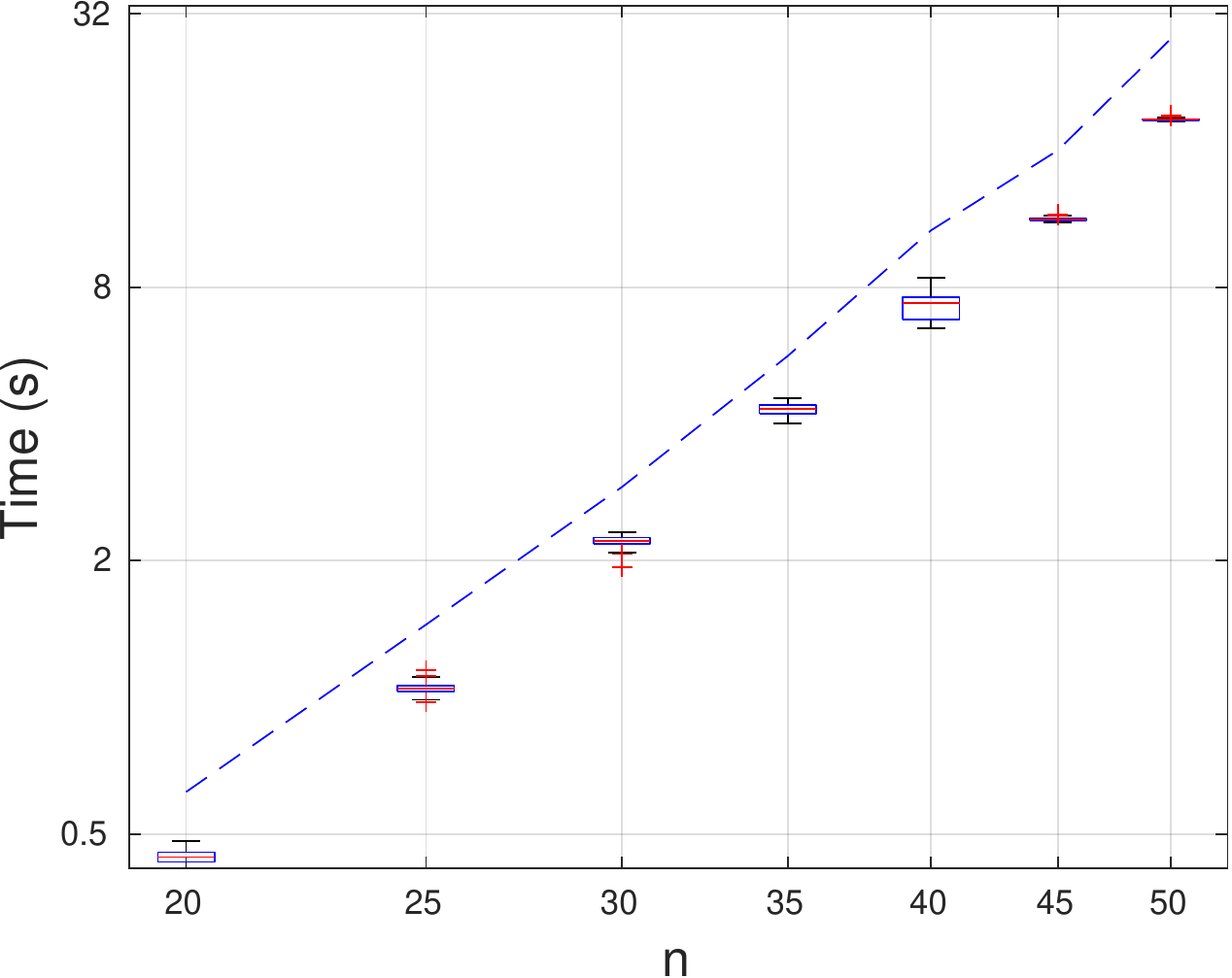}}
    \subfigure[Test errors on $\Xi_{n_i}$.]{\label{fig:3d}\includegraphics[width=0.327\linewidth]{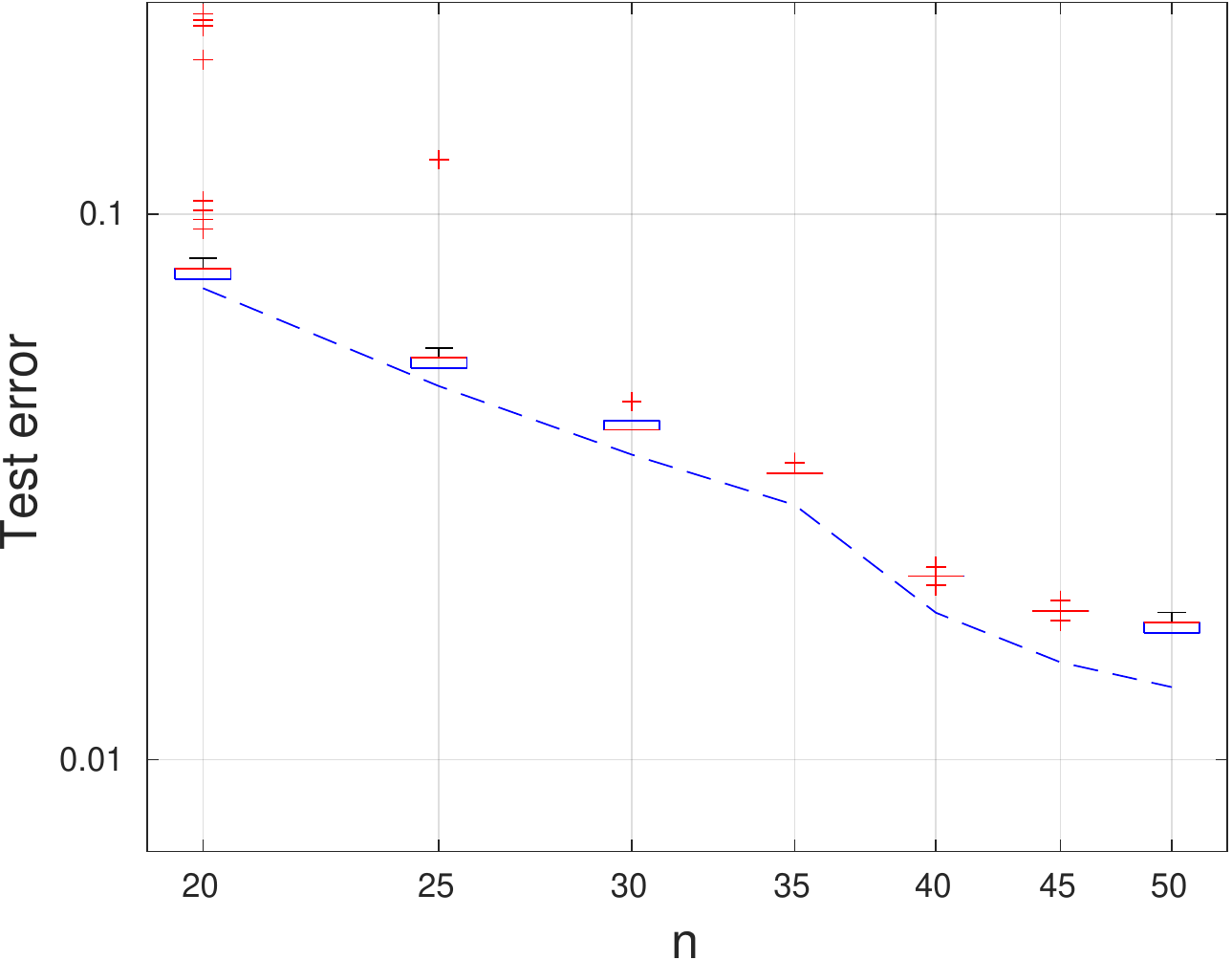}}
        \caption{Test 2. On the $x$-axis, we let $n_1,\dots,n_{\ell}$ vary. The results related to ERA are depicted by means of a dashed blue line, while the performance of SERA with $q=0.3$ are reported by using boxplots. We consider the kernel $\varphi_{W,\varepsilon}$ and the interpolation sets $H_{n_i}$, $i=1,\dots,\ell$. In Figures \ref{fig:3b} and \ref{fig:3d} both axes are in logarithmic scale.}
        \label{fig:3}
\end{figure}

The experiments provided in Section \ref{sec:testuno} suggest that the reduction ratio
$q=0.2$ might be taken into account as a reliable trade-off between computational time and accuracy. In order to investigate on this heuristic intuition in a different test, in the following we fix such a value for the reduction ratio, we set $p=3$, we take different Halton's point interpolation sets $H_{n_i}$ in $\Omega$, $i=1,\dots,\ell$, with $n_i$ equispaced between $n_1=20$ and $n_{\ell}=50$, $\ell=7$, and corresponding equispaced evaluation grids $\Xi_{n_i}$ in $\Omega$. Moreover, we consider a different interpolation task by taking the test function $g:\Omega\longrightarrow\mathbb{R}$ defined as
\begin{equation*}
    g(\bs{x})=x_1^2-x_2^4+e^{-(x_1+x_2)^2},\quad \bs{x}=(x_1,x_2),
\end{equation*}
and the \textit{Wendland $C^2$} kernel (see e.g. \cite[Section 11.2]{Fasshauer07}) $\varphi_{W,\varepsilon}(r)= (1-\varepsilon r)_+^4(4 \varepsilon  r +1)$
with $\bs{\varepsilon}_3 \in [0.01,2]$, uniformly, as vector of $101$ shape parameter values. Here, we use Tikhonov regularization with regularizing parameter $\lambda=10^{-10}$. In Figure \ref{fig:3}, the results confirm the suitability of the chosen reduction ratio according to Section \ref{sec:testuno}.
Numerical results for other values of 
$q\in[0.1,0.3]$ yield subplots with similar trends, which indicate 
that the tuning is not critical, and are omitted from this report.

\section{Discussion and conclusions}

In this paper, we proposed a stochastic approximation of the ERA that is based upon a low-rank approximation of the full kernel matrix inverse $\mathsf{K}_{\varepsilon}^{-1}$ required by the scheme. In Section \ref{sec:numerics}, we carried out some numerical experiments to compare the performance of the ERA and the proposed SERA, which has been detailed in Section \ref{sec:sera}. As confirmed by the tests, in which we considered different settings both regularized and non-regularized,
the proposed SERA in comparison to ERA
\begin{itemize}
  \item can select nearby shape parameters of the same magnitude, but
  \item with a saving in computational cost, especially as the rank of the approximating matrix $\mathsf{V}_s$ gets relatively small with respect to the one of $\mathsf{K}_{\varepsilon}^{-1}$, and
  \item can result in smaller interpolation error, i.e., better interpolants, with high probabilities.
\end{itemize}
Therefore, the SERA may be considered for fast shape parameter tuning in the context of RBF approximation. Future work consists of further investigations concerning the optimization of the reduction ratio parameter.

\section*{Acknowledgements}
This research has been accomplished with the financial support of GNCS-IN$\delta$AM and partially funded by the ASI-INAF grant \lq\lq Artificial Intelligence for the analysis of solar FLARES data (AI-FLARES)\rq\rq$\;$and the Hong Kong Research Grant Council GRF Grants.


\begin{thebibliography}{99}


\bibitem{Bekas07}
\textsc{C. Bekas, E. Kokiopoulou, Y. Saad}, \emph{An estimator for the diagonal of a matrix}, Appl. Num.
Math. \textbf{57}(11) (2007), pp. 1214--1229.

\bibitem{Cavoretto19}
\textsc{R. Cavoretto, A. De Rossi, M.S. Mukhametzhanov et al.}, \emph{On the search of the shape parameter in radial basis functions using univariate global optimization methods}, J. Glob. Optim. (2019).

\bibitem{Cavoretto18}
\textsc{R. Cavoretto, A. De Rossi, E. Perracchione}, \emph{Optimal selection of local approximants in RBF-PU interpolation}, J Sci Comput \textbf{74} (2018), pp. 1--22.

\bibitem{Celisse08}
\textsc{A. Celisse, S. Robin}, \emph{Nonparametric density estimation by exact leave-p-out cross-validation}, CSDA \textbf{52}(5) (2008), pp. 2350--2368.


\bibitem{Fasshauer07}
\textsc{G.E. Fasshauer},
\emph{Meshfree approximations methods with \textsc{Matlab}},
World Scientific, Singapore, 2007.

\bibitem{Fasshauer07a}
\textsc{G.E. Fasshauer, J.G. Zhang}, \emph{On choosing \lq\lq optimal\rq\rq  ~shape parameters for RBF approximation}, Numer. Algorithms \textbf{45} (2007), pp. 345--368.

\bibitem{Fasshauer15}
\textsc{G.E. Fasshauer, M.J. McCourt},
\emph{Kernel-based approximation methods using} \textsc{Matlab},
World Scientific, Singapore, 2015.

\bibitem{Fornberg07}
\textsc{B. Fornberg, J. Zuev}, \emph{The Runge phenomenon and spatially variable shape parameters in RBF interpolation}, Comput. Math. Appl. \textbf{54}(3) (2007), pp. 379--398.

\bibitem{Golub79}
\textsc{G. H. Golub, M. Heath, G. Wahba}, \emph{Generalized cross-validation as a method
for choosing a good ridge parameter}, Technometrics \textbf{21}(2) (1979), pp. 215--223.

\bibitem{Halton60}
\textsc{J. H. Halton}, \emph{On the efficiency of certain quasi-random sequences of points in evaluating multi-dimensional integrals}, Numer. Math. \textbf{2} (1960), pp. 84--90.

\bibitem{Ling16}
\textsc{L. Ling},
\emph{A fast block-greedy algorithm for quasi-optimal meshless trial subspace selection}, SIAM J.Sci. Comput. \textbf{38}(2) (2016), A1224–A1250.

\bibitem{Marchetti21}
\textsc{F. Marchetti},
\emph{The extension of Rippa’s algorithm beyond LOOCV}, Appl. Math. Lett. \textbf{120} (2021), 107262.

\bibitem{Mongillo11}
\textsc{M. Mongillo}, \emph{Choosing basis functions and shape parameters for radial basis function methods}, SIAM SIURO publications \textbf{4} (2011).

\bibitem{Rippa99}
\textsc{S. Rippa}, \emph{An algorithm for selecting a good value for the parameter $c$ in radial basis function interpolation}, Adv. Comput. Math. \textbf{11} (1999), pp. 193--210.

\bibitem{Scheuerer11}
\textsc{M. Scheuerer}, \emph{An alternative procedure for selecting a good value for the parameter c in RBF-interpolation}, Adv. Comput. Math. \textbf{34} (2011), pp. 105--126.

\bibitem{Yang18}
\textsc{F. Yang, L. Yan, L. Ling}, \emph{Doubly stochastic radial basis function methods}, J. Comput. Phys. \textbf{363} (2018), pp. 87--97.


\end{thebibliography}
\end{document}